# ON RIGIDITY AND THE ALBANESE VARIETY FOR PARALLELIZABLE MANIFOLDS

JÖRG WINKELMANN

ABSTRACT. We study the rigidity questions and the Albanese Variety for Complex Parallelizable Manifolds. Both are related to the study of the cohomology group $H^1(X, \mathcal{O})$. In particular we show that a compact complex parallelizable manifold is rigid iff $b_1(X) = 0$ iff $Alb(X) = \{e\}$ iff $H^1(X, \mathcal{O}) = 0$.

## 1. INTRODUCTION

For a compact Kähler manifold $X$ there exists a compact complex torus $Alb(X)$ and a holomorphic map $\pi : X \to Alb(X)$ such that

1. For every holomorphic map $f$ from $X$ to a compact complex torus $T$ there exists an element $a \in T$ and a holomorphic Lie group homomorphism $F : Alb(X) \to T$ such that $f(x) = a \cdot F \circ \pi(x)$ for all $x \in X$.
2. The holomorphic map $\pi : X \to Alb(X)$ induces isomorphisms on the cohomology groups $H^1(\cdot, \mathbb{R})$, $\Gamma(\cdot, d\mathcal{O})$ and $H^1(\cdot, \mathcal{O})$.

We prove similar statements for parallelizable manifolds. A complex manifold is called parallelizable if its tangent bundle is holomorphically trivial. A compact complex manifold is parallelizable if and only if it is biholomorphic to a quotient $G/\Gamma$ of a complex Lie group $G$ by a discrete subgroup $\Gamma$ ([21]). Compact complex parallelizable manifolds are never Kähler unless they are tori ([21]). Nevertheless, concerning the Albanese variety it is possible to obtain results similar to the Kähler case as long as reductive factors of rank 1 are absent from the Lie group $G$.

More precisely we prove the following.

**Theorem 1.** *Let $G$ be a connected complex linear-algebraic group and $\Gamma$ a lattice (i.e. a discrete subgroup such that the quotient manifold $X = G/\Gamma$ has finite volume with respect to Haar measure). Assume that there exists no surjective algebraic group homomorphism from $G$ to $PSL_2(\mathbb{C})$ or $\mathbb{C}^*$.*

*Then there exists a compact complex torus $A$ and a surjective holomorphic map $\pi : X \to A$ such that*

1. *For every connected complex Lie group $T$ and every holomorphic map $f : X \to T$ there exists an element $a \in T$ and a holomorphic Lie group homomorphism $F : A \to T$ such that $f(x) = a \cdot F(\pi(x))$ for all $x \in X$.*

Part of the work was done during the special year on Several Complex Variables at the MSRI. Research at MSRI is supported in part by NSF grant DMS-9022140. The author also wishes to thank the Harvard University for their invitation.

2. The map $\pi : X \to A$ induces isomorphisms on the cohomology groups $H^1(\cdot, \mathbb{C})$ and $\Gamma(\cdot, d\mathcal{O})$.
3. Every meromorphic function on $X$ is a pull-back from $A$.
4. If $X$ is compact, $\pi$ induces an isomorphism of $H^1(\cdot, \mathcal{O})$.

**Corollary 1.** *Let $G$ be a connected complex Lie group and $\Gamma$ a lattice. Then there exists an abelian variety $A$ and a surjective holomorphic connected $G$-equivariant map $\pi : G/\Gamma \to A$ inducing an isomorphism of the respective function fields.*

Furthermore, we determine $H^1(X, \mathcal{O})$ for cocompact lattices in arbitrary complex linear-algebraic groups.

**Theorem 2.** *Let $G$ be a connected complex linear-algebraic group, $\Gamma \subset G$ a lattice and $X = G/\Gamma$. Let $G = S \cdot R$ be a Levi-Malcev-decomposition, $N$ the nilradical, and $A = [S, R] \cdot N'$. In addition, let $W$ denote the maximal subvectorspace of $\mathcal{L}ie(R'A/A)$ such that $Ad(\gamma)|_W$ is a semisimple linear endomorphism with only real eigenvalues for every $\gamma \in \Gamma$.*

*Then $\dim H^1(G/\Gamma, \mathcal{O}) = \dim(G/G') + b_1(\Gamma/(R \cap \Gamma)) + \dim(W)$.*

Using this result we may deduce the following vanishing theorem.

**Theorem 3.** *Let $G$ be a connected complex Lie group, $\Gamma$ a discrete cocompact subgroup and $X = G/\Gamma$.*

*Then $H^1(X, \mathcal{O}) = \{0\}$ iff $b_1(X) = 0$.*

Since $H^1(X, \mathcal{T}) \simeq H^1(X, \mathcal{O}^n)$ for an $n$-dimensional parallelizable manifold, this implies that a compact complex parallelizable manifold admits infinitesimal deformations if and only if $b_1(X) > 0$. In addition, we will prove that for $b_1(X) > 0$ there actually exist small deformations. Thus we obtain the following result characterizing rigidity.

**Theorem 4.** *Let $G$ be a connected complex Lie group, $\Gamma$ a discrete cocompact subgroup and $X = G/\Gamma$ the complex quotient manifold. Then the following properties are equivalent.*

1. *$X$ admits no infinitesimal deformations (i.e. there are no non-trivial families over $Spec\, \mathbb{C}[\epsilon]/(\epsilon^2)$.*
2. *$X$ admits no small deformations (i.e. there are no non-trivial families over the unit disk).*
3. *$b_1(X) = 0$.*
4. *$H^1(X, \mathcal{O}) = \{0\}$.*

This characterization of rigidity of parallelizable manifolds generalizes a number of earlier results. The existence of non-trivial deformations of positive-dimensional tori is classical. Deformations of solv-manifolds have been studied by Nakamura [16]. Raghunathan proved the rigidity of quotients of $G/\Gamma$ for $G$ semisimple without *rank* 1-factors. Ghys proved theorem 4 for $G = SL_2(\mathbb{C})$.

**Theorem 5.** *Let $G$ and $H$ be linear-algebraic Let $G, H$ be a complex linear-algebraic groups, $\Gamma$ and $\Lambda$ discrete cocompact subgroups of $G$ resp. $H$ and $X = G/\Gamma$ and $Y = H/\Lambda$ the respective quotient manifolds.*

*Assume that* $\pi_1(X) \simeq \pi_1(Y)$.
*Then* $\dim H^1(X, \mathcal{O}) = \dim H^1(Y, \mathcal{O})$.

The assumption of $G$ being linear-algebraic is not very restrictive (in fact we do not know any compact complex parallelizable manifold for which it can not be fulfilled), but it is crucial that both $X$ and $Y$ are parallelizable. Nakamura gave an example where $X$ is a compact complex parallelizable manifold and $Y$ a non-parallelizable compact complex manifold diffeomorphic to $X$ (thus in particular $\pi_1(X) \simeq \pi_1(Y)$), but $\dim H^1(X, \mathcal{O}) \neq \dim H^1(Y, \mathcal{O})$ [16].

Parts of the results presented here are already contained in the Habilitationsschrift of the author [22].

## 2. Preparations

We summarize some basic results which are used in the sequel.

**Theorem 6.** *(Malcev, [13]) Let $G$ be a real nilpotent Lie group, $\Gamma$ a discrete subgroup such that $G/\Gamma$ has finite volume.*
*Then $G/\Gamma$ is compact, $\Gamma'$ is a lattice in $G'$ and $\Gamma'$ is of finite index in $G' \cap \Gamma$.*

**Theorem 7.** *(Kodaira, [16]) Let $G$ be a complex nilpotent Lie group and $\Gamma$ a lattice. Then $H^1(G/\Gamma, \mathcal{O}) \simeq H^1(G/G'\Gamma, \mathcal{O})$.*

**Theorem 8.** *(Margulis, [14]) Let $S$ be a semisimple complex Lie group without rank 1-factor and $\Gamma$ a lattice. Then $\Gamma'$ is of finite index in $\Gamma$.*

**Theorem 9.** *(Mostow, [15]) Let $G$ be a complex connected Lie group, $R$ its radical and $N$ its nilradical and $\Gamma$ a lattice.*
*Then $G\Gamma$ and $N\Gamma$ are closed and $R/(R \cap \Gamma)$ and $N/(N \cap \Gamma)$ are compact.*

The following theorem of Akhiezer generalizes earlier results of Raghunathan.

**Theorem 10.** *(Akhiezer, [2]) Let $G$ be a simply-connected semisimple complex Lie group and $\Gamma$ a cocompact lattice. Let $X = G/\Gamma$.*
*In this case the natural morphism $i : H^1(X, \mathbb{C}) \to H^1(X, \mathcal{O})$ is an isomorphism.*

**Proposition 1.** *Let $G$ be a connected complex Lie group, $\Gamma$ a lattice and $H$ a closed complex Lie subgroup of $G$ containing $\Gamma$.*
*Then the connected component $H^0$ of $H$ is normal in $G$ and therefore $G/H$ is again parallelizable.*

*Proof.* We consider the Tits-fibration, i.e. the natural $G$-equivariant morphism from $G/H$ to a Grassmann manifold given by $x \mapsto \mathcal{L}ie(G_x)$. From Iwamoto's generalization of the Borel density theorem (see [11]) it follows that this map is constant. This implies that $H^0$ is normal. □

## 3. Orbits of the commutator group

**Lemma 1.** *Let $G$ be a complex Lie group, $\Gamma$ a discrete subgroup and $A$ a normal simply-connected commutative complex Lie subgroup of $G$. Assume that $A/(A \cap \Gamma)$ is compact and that the Zariski-closures of $Ad(\Gamma)$ and $Ad(G)$ in $GL(\mathcal{L}ie\, G)$ coincide.*
*Then $[\Gamma, A \cap \Gamma]$ is a lattice of $[G, A]$ and a subgroup of finite index in $[G, A] \cap \Gamma$.*

*Proof.* The group $A$ is a complex vector space and conjugation in the group $G$ induces a representation $\rho : G/A \to GL(A)$. This representation is naturally isomorphic to the restriction of the adjoint representation $Ad$ to $\mathcal{L}ie\, A$.

For $g \in G$ define $\zeta_g \in End(A)$ by $\zeta_g = \rho(g) - id_A$. Let $\Lambda = \oplus_{\gamma \in \Gamma} \zeta_\gamma (A \cap \Gamma)$. Clearly, $[\Gamma, A \cap \Gamma]$ is the subgroup generated by $\Lambda$. Since $A \cap \Gamma$ is cocompact in $A$, $<A \cap \Gamma>_\mathbb{R} = <A \cap \Gamma>_\mathbb{C}$. This implies that $<\Lambda>_\mathbb{R} = <\Lambda>_\mathbb{C}$ and consequently that $<\Lambda>_\mathbb{Z}$ is a lattice in the complex subspace $V = <\Lambda>_\mathbb{C}$. We claim that $[A, G'] = V$. Indeed, $Ad(\gamma)(\mathcal{L}ie\, A) \subset \mathcal{L}ie\, V$ for all $\gamma \in \Gamma$ implies that $Ad(g)(\mathcal{L}ie\, A) \subset \mathcal{L}ie\, V$ for all $g \in G$, because we assumed that the Zariski-closures of $Ad(\Gamma)$ and $Ad(G)$ coincide. □

**Proposition 2.** *Let $G$ be a simply-connected complex Lie group, $R$ its radical and $\Gamma$ a lattice.*

*Then $(G' \cap R)\Gamma$ is closed in $G$, $R \cap \Gamma'$ is a subgroup of finite index in $G' \cap R \cap \Gamma$ and a cocompact lattice in $G' \cap R$.*

*Proof.* First recall that a normal Lie subgroup in a simply-connected Lie group is also simply-connected.

Let $N$ denote the nilradical and $N^k$ the central series of $N$. By the theorems of Mostow and Malcev all of the groups $N$ and $N^k$ ($k \in \mathbb{N}$) have compact orbits in $G/\Gamma$. Furthermore $N^k/N^{k+1}$ is a simply-connected commutative normal Lie subgroup of $G/N^{k+1}$ for all $k$.

Thus we may apply the preceding lemma repeatedly and obtain that $G' \cap N \cap \Gamma'$ is of finite index in $G' \cap N \cap \Gamma$ and a lattice in $N \cap G'$.

This is the statement of the proposition, since $G' \cap R \subset N$. □

**Corollary 2.** *(Barth-Otte [5]) Let $G$ be a solvable complex Lie group and $\Gamma$ a lattice. Then $G'\Gamma$ is closed in $G$.*

**Theorem 11.** *Let $G$ be a simply-connected complex Lie group, $R$ its radical and $\Gamma$ a lattice. Assume that $\Gamma'/(R \cap \Gamma')$ is a subgroup of finite index in $\Gamma/(R \cap \Gamma)$.*

*Then $G'\Gamma$ is closed in $G$, $G' \cap \Gamma$ is a lattice in $G'$ and $\Gamma'$ is a subgroup of finite index in $G' \cap \Gamma$.*

**Corollary 3.** *Let $G$ be a simply-connected complex Lie group, $R$ its radical and $\Gamma$ a lattice. Assume that no simple factor of $G/R$ is isomorphic to $SL_2(\mathbb{C})$.*

*Then $G'\Gamma$ is closed in $G$, $G' \cap \Gamma$ is a lattice in $G'$ and $\Gamma'$ is a subgroup of finite index in $G' \cap \Gamma$.*

Recursive application of the theorem implies the following.

**Corollary 4.** *Let $G$ and $\Gamma$ be as in the above theorem. Let $G^k$ denote the $k$-th term of the derived series of $G$. Then $G^k\Gamma$ is closed in $G$, $G^k \cap \Gamma$ is a lattice in $G^k$ and $\Gamma^k$ is a subgroup of finite index in $G^k \cap \Gamma$.*

*Proof.* By the above proposition it suffices to consider the case where $R \cap G' = \{e\}$. Then $R$ is central in $G$. Since $G$ is simply-connected, this implies $G = R \times S$ for a semisimple complex Lie group $S$. Let $\pi : G \to S$ denote the natural projection. Then $\pi(\Gamma)$ is a lattice in $S$. Margulis' theorem implies that $\pi(\Gamma')$ is also a lattice in $S$. It follows that $(\Gamma \cap R) \cdot (\Gamma \cap S) \supset (\Gamma \cap R) \cdot \Gamma'$ is already a lattice in $G$. In particular, $\Gamma'$ is a lattice in $S$ and hence of finite index in $S \cap \Gamma = G' \cap \Gamma$. □

Although in general $G' \cap \Gamma$ is not a lattice in $G'$, it is still large in a certain sense.

**Proposition 3.** *Let $G$ be a complex Lie group, $\Gamma$ a lattice in $G$.*
*Then every $\Gamma'$-invariant plurisubharmonic function on $G'$ is constant.*

*Proof.* Since $(R \cap G')/(R \cap \Gamma')$ is compact, it is clear that $\Gamma'$-invariant plurisubharmonic functions are constant along the $R \cap G'$-orbits. Thus we may assume that $R \cap G' = \{e\}$. Then $G$ is a direct product of an abelian complex Lie group $R$ and a semisimple complex Lie group $S$. Recall that a lattice is Zariski-dense. This implies that $\Gamma'$ is Zariski-dense in $G'$, because for sufficiently large $n$ the map $\zeta : (G \times G)^n \to G'$ given by

$$\zeta : ((x_1, y_1), \ldots, (x_n, y_n)) \mapsto (x_1 y_1 x_1^{-1} y_1^{-1}) \cdot \ldots \cdot (x_n y_n x_n^{-1} y_n^{-1}) \tag{1}$$

is a dominant regular morphism and therefore maps $(\Gamma \times \Gamma)^n$ onto a Zariski-dense subset of $G'$. Now a theorem of Berteloot and Oeljeklaus [4] implies that $\Gamma'$-invariant plurisubharmonic functions must be $S$-invariant. □

In a similar way one can prove the following.

**Proposition 4.** *Let $G$ be a connected complex Lie group and $\Gamma$ a lattice.*
*Then every plurisubharmonic function on $X = G/\Gamma$ is constant.*

## 4. The Construction of the Albanese torus

**Proposition 5.** *Let $G$ be a complex Lie group and $\Gamma$ a closed complex Lie subgroup. Assume that every holomorphic function on $G/\Gamma$ is constant.*

*Let $I$ denote the smallest closed complex Lie subgroup of $G$ containing $G'\Gamma$.*

*Then $G/I$ is a commutative complex Lie group and the natural projection $\pi : G/\Gamma \to G/I$ has the following properties.*

1. *For every connected complex Lie group $T$ and every holomorphic map $f : G/\Gamma \to T$ there exists an element $a \in T$ and a holomorphic Lie group homomorphism $F : G/I \to T$ such that $f(x) = a \cdot F \circ \pi(x)$ for all $x \in X$.*
2. *If $I = G'\Gamma$, then every closed holomorphic 1-form on $G/\Gamma$ is a pull-back from $G/I$.*
3. *If $\Gamma$ is a lattice in $G$, then $G/I$ is compact.*
4. *If $\Gamma$ is a lattice in $G$, then every meromorphic function on $G/\Gamma$ is a pull-back from $G/I$.*

*Proof.* 1. Let $Ad$ denote the adjoint representation of $T$ and $Z$ the connected component of the center of $T$. Since holomorphic functions on $GL(\mathcal{L}ie\,T)$ separate points, it is clear that $Ad \circ f$ is constant for every holomorphic map $f : G/\Gamma \to T$. Thus there is no loss in generality in assuming $Z = T$.

The cotangent bundle $T^*(Z)$ is spanned by invariant closed holomorphic 1-forms $\omega_i$. Now

$$0 = df^*\omega_i(X, Y) = f^*\omega_i([X, Y]) + X(f^*\omega_i)Y - Y(f^*\omega_i X) = f^*\omega_i([X, Y])$$

for every holomorphic map $f: G/\Gamma \to T$, every closed holomorphic 1-form $\omega_i$ and every holomorphic vector fields $X, Y$ (Note that $(f^*\omega_i Y)$ is a holomorphic function on $G/\Gamma$ and therefore constant.)

It follows that the $G'$-orbits are contained in the fibers of $f$.

The map $f$ must be equivariant [23]. Hence $f$ fibers through $G/I$.

2. The equation $(*)$ implies that $\Gamma(X, d\mathcal{O}) \simeq \mathcal{L}ie(G/G')^*$. Hence the assertion.

3. By a Fubini-type argument it follows that $G/I$ has finite volume. (see [20], Lemma 1.6). Since $G/I$ is an abelian group, this implies compactness.

4. As a homogeneous complex manifold, $G/\Gamma$ admits a meromorphic reduction $G/\Gamma \to G/J$ with meromorphically separable base $G/J$ ([10]). The base manifold $G/J$ is again parallelizable (prop. 1). Therefore it suffices to show the following:

**Lemma 2.** *Let $H$ be a connected complex Lie group and $\Lambda$ a lattice. Assume that $Y = H/\Lambda$ is meromorphically separable.*

*Then $H$ is commutative and $Y$ a compact complex torus.*

*Proof.* Recall that the radical $R$ of $H$ has closed orbits in $Y = H/\Lambda$ and that these orbits are moreover compact, because lattices in solvable Lie groups are necessarily cocompact. Thus we obtain a holomorphic fiber bundle $\tau: Y \to Z = H/R\Lambda$ with compact fibers. Using meromorphic separability of $Y$, for any point $y \in Y$ and any natural number $k$ with $0 \le k \le \dim Y$ we can find a closed analytic subset $W \subset Y$ of pure dimension $k$ containing $y$ (simply by taking irreducible components of intersections of zero divisors of meromorphic functions on $Y$). Since the projection map $\tau: Y \to Z$ is a proper holomorphic map, this implies the existence of analytic hypersurfaces in $W$ unless $W$ is trivial (i.e. a point).

Recall that a theorem of Huckleberry and Margulis [9] states a quotient of a complex semisimple Lie group by a Zariski-dense subgroup never contains analytic hypersurfaces. Hence in our situation $W$ must be a point, i.e. $H$ is solvable. It follows that $H/\Lambda$ is compact. Finally observe that a meromorphically separable compact homogeneous manifold is automatically projective and in particular Kähler ([8]). But a parallelizable compact complex manifold is Kähler if and only if it is a compact complex torus ([21]). □

□

## 5. Leray spectral sequence

Let $f: X \to Y$ be a holomorphic map between complex spaces. There is a Leray spectral sequence for the sheaf $\mathcal{O}_X$. The respective lower term sequence yields the following.

$$0 \to H^1(Y, \mathcal{R}^0 f_* \mathcal{O}_X) \to H^1(X, \mathcal{O}_X) \to H^0(Y, \mathcal{R}^1 f_* \mathcal{O}_X) \to H^2(Y, \mathcal{R}^0 f_* \mathcal{O}_X)$$

Assume that $f$ is connected and proper. Then $\mathcal{R}^0 f_* \mathcal{O}_X = \mathcal{O}_Y$. Furthermore, if $\dim H^k(f^{-1}\{p\}, \mathcal{O}) = r$ for a natural number $k$ and all $p \in Y$, then, by Grauert's theorem, $\mathcal{R}^k f_* \mathcal{O}_X$ is a locally free coherent sheaf of rank $r$. This implies in particular the following observation.

**Lemma 3.** *Let $f : X \to Y$ be a proper connected holomorphic map and assume that for every fiber $F_p$ the induced cohomology map $i^* : H^1(X, \mathcal{O}) \to H^1(F_p, \mathcal{O})$ vanishes. Assume furthermore that $\dim H^1(F_p, \mathcal{O})$ does not depend on $p \in Y$.*

*Then $\alpha : H^1(X, \mathcal{O}) \to H^0(\mathcal{R}^1 f_* \mathcal{O})$ is the zero map.*

## 5.1. Leray spectral sequence for flat bundles.

**Proposition 6.** *Let $\pi : E \to X$ be a holomorphic fiber bundle with typical fiber $F$ and structure group $S$. Assume that $F$ is connected and compact, that $S$ acts on $F$ in such a way that there exists an invariant hermitian metric and that $E \to X$ admits a flat holomorphic connection.*

*Then $\alpha : H^1(E, \mathcal{O}) \to H^0(X, \mathcal{R}^1 \pi_* \mathcal{O})$ is surjective.*

*Proof.* Let $\eta \in H^0(X, \mathcal{R}^1 \pi_* \mathcal{O})$. Then there exists an open cover $\mathcal{U} = (U_i)_i$ by contractible open Stein subsets of $X$ such that $\eta$ is given by $\eta_i \in H^1(\pi^{-1}(U_i), \mathcal{O})$. Using the Dolbeault-isomorphism we may choose corresponding $\bar\partial$-closed $(0,1)$-forms $\omega_i$ on $W_i = \pi^{-1}(U_i)$. Of course these forms $\omega_i$ are not unique. Now, if there is a *canonical* way to choose $\omega_i$, then the forms $\omega_i$ coincide on the intersection of the $U_i$ and yield a globally defined $\bar\partial$-closed $(0,1)$-form on $E$, implying that $\eta$ is contained in the image of $H^1(E, \mathcal{O})$.

The assumptions made in the proposition allow us to choose a hermitian metric on each fiber in such a way that for every contractible Stein open subset $U \subset X$ we obtain a trivialization $E|_U \simeq U \times F$ which is compatible with the flat connection and such that the choosen hermitian metric on each fiber is just the pull-back of one fixed $S$-invariant hermitian metric on $F$. Then there is a canonical way to choose the forms $\omega_i$. Namely the forms $\omega_i$ are to be choosen such that they annihilate horizontal vector fields (with respect to the connection) and are harmonic if restricted to a fiber of $\pi$. $\square$

**Lemma 4.** *Let $G$ be a complex Lie group, $\Gamma$ a discrete subgroup and $A$ a connected commutative normal complex Lie subgroup such that $A/(A \cap \Gamma)$ is compact. Assume that the short exact sequence of Lie algebras*
$$0 \to \mathcal{L}ie(A) \to \mathcal{L}ie(G) \to \mathcal{L}ie(G/A) \to 0$$
*is split.*

*Then $\pi : E = G/\Gamma \to G/A\Gamma = B$ is a torus bundle with flat holomorphic connection.*

*Proof.* The flat connection is induced by the splitting of the Lie algebra sequence. $\square$

**Corollary 5.** *Under the assumptions of the lemma the induced map $H^1(E, \mathcal{O}) \to H^0(B, \mathcal{R}^1 \pi_* \mathcal{O})$ is surjective.*

We will also need a description of the structure of $\mathcal{R}^1 f_* \mathcal{O}$ for torus bundles.

Let $f : E \to Y$ be a locally trivial holomorphic fiber bundle with a compact complex torus $T$ is typical fiber. Let $V = \Omega^1(T)$ denote the vector space of holomorphic 1-forms on $T$. Then there is a exact sequence
$$1 \to \underbrace{Aut^0(T)}_{=T} \to Aut(T) \xrightarrow{\zeta} GL(V)$$

Let $\mathcal{U} = \{U_i\}$ be a trivializing open cover of $Y$ such that $E$ is given by transition functions $\phi_{ij} : U_i \cap U_j \to Aut(T)$.

**Claim 1.** *Under the above assumptions $\mathcal{R}^1 f_* \mathcal{O}_E$ is a locally free coherent sheaf on $Y$ isomorphic to the sheaf of sections in the flat vector bundle given by the transition functions $\psi_{ij} : U_i \cap U_j \to GL(V)$ defined by $\psi_{ij} = \overline{\zeta \circ \phi_{ij}}$.*

*Proof.* This is a consequence of the Dolbeault-isomorphism. □

We apply this to parallelizable manifolds.

**Proposition 7.** *Let $G$ be a complex Lie group, $\Gamma$ a discrete subgroup and $A$ a normal abelian complex Lie subgroup. Assume that $A/(A \cap \Gamma)$ is compact. Denote the natural projection map $E = G/\Gamma \to G/A\Gamma = B$ by $\pi$.*

*In this case $\mathcal{R}^1 \pi_* \mathcal{O}$ is a flat vector bundle of rank $\dim A$ over $B$ induced by the representation $\rho : \Gamma \to GL(\mathcal{L}ie\, A^*)$ given by $\gamma \mapsto \overline{Ad^*(\gamma)}$.*

**Lemma 5.** *Let $\pi : X \to Y$ be a finite holomorphic covering.*
*Then $\pi^* : H^1(Y, \mathcal{O}) \to H^1(X, \mathcal{O})$ is injective.*

*Proof.* If $\omega$ is a $\bar{\partial}$-closed $(0,1)$-form on $Y$ and $\bar{\partial} f = \pi^* \omega$ for a function $f$ on $X$, then $\bar{\partial} g = \omega$ for $g(y) = \frac{1}{d} \sum_{\pi(x)=y} f(x)$ (where $d$ denotes the degree of $\pi$). □

## 6. Description of $H^1(X, \mathcal{O})$

**Theorem 12.** *Let $G$ be a connected complex Lie group, $\Gamma \subset G$ a lattice and $X = G/\Gamma$. Let $G = S \cdot R$ be a Levi-Malcev-decomposition, $N$ the nilradical, and $A = [S, R] \cdot N'$. Furthermore let $W$ denote the maximal linear subspace of $\mathcal{L}ie(R'A/A)$ such that $Ad(\gamma)|_W$ is a semisimple linear endomorphism with only real eigenvalues for every $\gamma \in \Gamma$.*

*Then $\dim H^1(G/\Gamma, \mathcal{O}) \leq \dim(G/G') + b_1(\Gamma/(R \cap \Gamma)) + \dim(W)$. Equality holds, if $G$ is linear-algebraic.*

**Remark 1.** *We do not know any example of a quotient manifold of a connected complex Lie group by a lattice which is not biholomorphic to a quotient of a linear-algebraic $\mathbb{C}$-group by a lattice.*

**Corollary 6.** *Let $G$ be a simply-connected complex Lie group and $\Gamma \subset G$ be a cocompact lattice. Assume that the radical $R$ is nilpotent and that the semisimple group $G/R$ contains no factor $S_0$ such that $S_0/(S_0 \cap R\Gamma)$ is compact and $S_0 \simeq SL_2(\mathbb{C})$.*
*Then $H^1(G/\Gamma, \mathcal{O}) \simeq H^1(G/G'\Gamma, \mathcal{O})$.*

*Proof.* If $R$ is nilpotent, then $R = N$ and therefore $A = R'A$. Since $A \subset H \subset R'A$, it follows that $A = H$. Furthermore the assumption on $S$ implies that $H^1(G/R\Gamma, \mathcal{O}) = 0$. Thus $\dim H^1(G/\Gamma, \mathcal{O}) = \dim G/G'$. The assumption on $S$ also implies that $G'\Gamma$ is closed in $G$. Thus $H^1(G/G'\Gamma, \mathcal{O}) \simeq H^1(G/\Gamma, \mathcal{O})$. □

We now prove the theorem.

*Proof.* We study the sequence of fibrations
$$1 \to G/\Gamma \xrightarrow{\pi_1} G/N'\Gamma \xrightarrow{\pi_2} G/A\Gamma \xrightarrow{\pi_3} G/(G' \cap R)\Gamma \xrightarrow{\pi_4} G/R\Gamma \to 1$$
First we must verify the existence of these fibrations, i.e. to show that $N'$, $A$, $G' \cap R$ and $R$ all have closed orbits. Results of Mostow [15] imply that $R$ and $N$ have closed orbits. By classical results of Malcev [13] this implies closedness of the $N'$-orbits. From prop. 2 we obtain closedness of the $G' \cap R$-orbits. Finally note that $A = (G^k \cap R)N'$ for $k$ sufficiently large. Hence closedness of $A$-orbits follows from thm. 11.

We note that all these projections $\pi_i$ are surjective proper holomorphic maps and for each $i$ we will study the lower term sequence of the Leray spectral sequence for the sheaf $\mathcal{O}$.

**Claim 2.** *The induced map* $\pi_1^*: H^1(G/N'\Gamma, \mathcal{O}) \to H^1(G/\Gamma, \mathcal{O})$ *is an isomorphism.*

*Proof.* By a theorem of Kodaira (see [16]) there is an isomorphism
$$H^1(N/N'(N \cap \Gamma), \mathcal{O}) \simeq H^1(N/(N \cap \Gamma, \mathcal{O}).$$
It follows that the embedding $i: N'/(N' \cap \Gamma) \hookrightarrow N/(N \cap \Gamma)$ induces the zero map between the respective cohomology groups $H^1(\cdot, \mathcal{O})$. Thus the group homomorphism $j^*: H^1(G/\Gamma, \mathcal{O}) \to H^1(N'/(N' \cap \Gamma))$ induced by the inclusion map $j: N'/(N' \cap \Gamma) \to G/\Gamma$ must vanish as well. By homogeneity it follows that for every fiber $F$ of $p$ the induced cohomology map $i^*: H^1(G/\Gamma, \mathcal{O}) \to H^1(F, \mathcal{O})$ is zero. Therefore
$$H^1(G/\Gamma, \mathcal{O}) \to H^0(G/N'\Gamma, \mathcal{R}^1 p_* \mathcal{O})$$
vanishes (lemma 3) and $H^1(G/\Gamma, \mathcal{O}) \simeq H^1(G/N'\Gamma, \mathcal{O})$. □

**Claim 3.** *The induced map* $\pi_2^*: H^1(G/A\Gamma, \mathcal{O}) \to H^1(G/N'\Gamma, \mathcal{O})$ *is an isomorphism.*

*Proof.* Recall that $A$ is defined as $A = [S, R]N'$. Since $[S, R] \subset R \cap G' \subset N$, it is clear that $A/N'$ is abelian. This implies that for $\pi: G/N'\Gamma \to G/A\Gamma$ the higher direct image sheaf $\mathcal{R}^1 \pi_* \mathcal{O}$ is the coherent sheaf associated to a flat vector bundle given by a representation $\rho$ which arises in the following way: $\rho$ is the restriction of the complex conjugation of the representation from $G$ on $GL(\mathcal{L}ie(A/N')^*)$ induced by the coadjoint representation. By construction no linear subspace of $\mathcal{L}ie(A/N')$ is invariant under $Ad(S)$. It follows that $\mathcal{R}^1 \pi_* \mathcal{O}$ does not admit global sections. Hence the claim. □

The lower term sequence of the Leray spectral sequence for $\pi_3$ yields
$$H^1(G/(G' \cap R)\Gamma, \mathcal{O}) \to H^1(G/A\Gamma, \mathcal{O}) \xrightarrow{\alpha} H^0(G/(G' \cap R)\Gamma, \mathcal{R}^1(\pi_3)_* \mathcal{O}) \xrightarrow{\phi} W$$
The isomorphism $\phi$ is a consequence of [24].

**Claim 4.** *If $G$ is linear-algebraic, then $\alpha$ is surjective.*

*Proof.* Using cor. 5 it suffices to show that the short exact sequence of complex Lie groups
$$(*) \qquad 1 \to (G' \cap R)/A \to G/A \to G/(G' \cap R) \to 1$$

is split. Let $V$ denote the unipotent radical of $G$. Then $G' \cap R \subset V \subset N$. Furthermore $A = [S,R]N'$ implies that $V/A$ is commutative. Let $L$ denote a Levi subgroup of $G$, i.e. a maximal connected reductive subgroup of $G$. Then $L$ is a reductive group, acting (by conjugation) linearly on the vector group $V/A$. Hence there is a $L$-invariant subvector space $W \subset V$ such that $V = W \oplus (G' \cap R)$. Thus $G = (L \ltimes W) \ltimes V$ and we obtain a splitting of $(*)$. □

Next we discuss the lower term sequence for $\pi_4$.

$$1 \to H^1(G/R\Gamma, \mathcal{O}) \to H^1(G/(G' \cap R)\Gamma, \mathcal{O}) \xrightarrow{\beta} H^0(G/R\Gamma, \mathcal{R}^1(\pi_4)_*\mathcal{O})$$

First we note that the exact sequence of Lie algebras

$$0 \to \mathcal{L}ie(R/(G' \cap R)) \to \mathcal{L}ie(G/(G' \cap R)) \to \mathcal{L}ie(G/R) \to 0$$

is always split. Hence $\beta$ is surjective. Furthermore the adjoint action of $G$ on $\mathcal{L}ie(G)$ induces the trivial action on $Lie(R/(G'\cap R))$. Thus $\mathcal{R}^1(\pi_4)_*\mathcal{O}$ is a free $\mathcal{O}$-module sheaf with rank equal to $\dim_{\mathbb{C}}(H^1(R/(G'\cap R)\Gamma, \mathcal{O})$. Note that $G/G' \simeq R/(G' \cap R)$ and that $R/(G' \cap R)$ is a compact complex torus. It follows that $\dim H^0(G/R\Gamma, \mathcal{R}^1(\pi_4)_*\mathcal{O}) = \dim(G/G')$.

Finally we recall that $\dim H^1(G/R\Gamma, \mathcal{O}) = \dim H^1(G/R\Gamma, \mathbb{C}) = b_1(\Gamma/(R \cap \Gamma))$ by a result of Akhiezer [2]. □

We now present a precise vanishing criterion for $H^1(G/\Gamma, \mathcal{O})$.

**Proposition 8.** *Let $G$ be a complex Lie group and $\Gamma$ a discrete cocompact subgroup. Then $H^1(G/\Gamma, \mathcal{O}) = 0$ if and only if $b_1(G/\Gamma) = 0$.*

*Proof.* The exponential sequence yields an embedding $H^1(X, \mathbb{Z}) \hookrightarrow H^1(X, \mathcal{O})$. Hence $H^1(X, \mathcal{O}) = 0$ implies $b_1(X) = 0$ for $X = G/\Gamma$.

Conversely let us assume $b_1(X) = 0$. Then $Hom(\Gamma, \mathbb{Z}) = 0$. Since $\Gamma$ is a lattice, this implies $Hom(G, \mathbb{C}) = 0$. Hence $G = G'$ and consequently $R = [S, R]$. This implies in particular that $R'A = A$ for $A = [S,R]N'$. Hence $b_1(\Gamma/(R\cap\Gamma))$, $\dim(G/G')$ and $\dim R'A/A$ all equal zero and $H^1(X, \mathcal{O}) = \{0\}$ follows by the above theorem. □

**Corollary 7.** *Let $X$ be a compact complex parallelizable manifold. Then $H^1(X, \mathcal{T}) = 0$ iff $b_1(X) = 0$.*

## 7. Topological invariance of $\dim H^1(X, \mathcal{O})$

For compact Kähler manifolds $\dim H^1(X, \mathcal{O})$ equals $b_1(X)$ and therefore depends only on the fundamental group of $X$. Nakamura has shown (thus answering a question of Iitaka) that for arbitrary compact complex manifolds $\dim H^1(X, \mathcal{O})$ may jump within a smooth family. He gave an example of a parallelizable complex manifold $X_0$ which admits small deformations $X_t$ for which $\dim H^1(X_t, \mathcal{O}) \ne \dim H^1(X, \mathcal{O})$. In this example the $X_t$ are no longer parallelizable for $t \ne 0$. This is not by coincidence. As an application of the description of $H^1(X, \mathcal{O})$ obtained in the preceding section we deduce the following result concerning the topological invariance of $H^1(X, \mathcal{O})$.

**Theorem 13.** *Let $G$, $H$ be a complex linear-algebraic groups, $\Gamma$ and $\Lambda$ discrete cocompact subgroups of $G$ resp. $H$ and $X = G/\Gamma$ and $Y = H/\Lambda$ the respective quotient manifolds.*

*Assume that $\pi_1(X) \simeq \pi_1(Y)$.*

*Then $\dim H^1(X, \mathcal{O}) = \dim H^1(Y, \mathcal{O})$.*

*Proof.* Let $G_0$ denote a complex-linear algebraic group, $G$ ist universal cover, $\Gamma_0$ a discrete cocompact subgroup of $G_0$ and $\Gamma$ its preimage in $G$ under the natural projection $G \to G_0$.

We have to show that $\dim H^1(G/\Gamma, \mathcal{O})$ is completely determined by $\Gamma$.

We observe that simply-connected complex Lie groups are linear and discrete subgroups in linear groups admit maximal normal solvable and maximal normal nilpotent subgroups which we denote by $rad(\cdot)$ resp. $n(\cdot)$. The usual density theorems for lattices in complex Lie groups imply that $rad(\Gamma)$ and $n(\Gamma)$ contain $R \cap \Gamma$ resp. $N \cap \Gamma$ as subgroups of finite index (where $R$ resp. $N$ denotes the radical resp. nilradical of $G$).

Among all subgroups of finite index of $n(\Gamma)$ we choose one for which the commutator group, which we call $n_1(\Gamma)$ has minimal $\mathbb{Z}$-rank (see e.g. [18], Def.2.9. for the notion of $\mathbb{Z}$-rank for finitely generated nilpotent groups). Then $n_1(\Gamma)$ is commensurable to $N' \cap \Gamma$. Similarly we choose $r_1(\gamma)$. Next recall that for a lattice $\Lambda$ in a semisimple Lie group $S$ the image $\rho(\Lambda)$ is Zariski-dense in $\rho(S)$. This holds for every *real* representation $\rho$. Besides $[S, N]N' = [G', N]N'$, since $G' \cap R \subset N$. It follows that $\Gamma_1 := [\Gamma', n(\Gamma)]n_1(\Gamma)$ must be a lattice in $[S, N]N'$.

Now $b_1(G/R\Gamma) = b_1(\Gamma/rad(\Gamma))$,

$$\dim G/G' = \dim R/(R \cap G') = \dim R/[S, N]R' = rank_{\mathbb{Z}}(rad(\Gamma)/\Gamma_1 r_1(\Gamma))$$

and $\dim_{\mathbb{C}}(W)$ equals half the dimension over $\mathbb{Q}$ of the $\mathbb{Q}$-vector space $W_0$ where $W_0$ is the maximal $\mathbb{Q}$-linear subspace of $(r_1(\Gamma)\Gamma_1/\Gamma_1) \otimes_{\mathbb{Z}} \mathbb{Q}$ such that the natural linear transformation $\tau(\gamma)$ induced by conjugation is diagonalizable over $\mathbb{R}$ with only real eigenvalues for every $\gamma \in \Gamma$.

Thus $\dim H^1(G/\Gamma, \mathcal{O})$ depends only on $\Gamma$. □

We want to emphasize that $\dim H^1(X, \mathcal{O})$ depends on $\pi_1(G/\Gamma)$ in a rather subtle way. As we will see (see section 8) the dimension of $H^1(X, \mathcal{O})$ may even jump within a given commensurability class.

## 8. Small deformations

So far we proved that there exist infinitesimal deformations iff $b_1(X) > 0$. Now we shall discuss small deformations. By actually constructing a deformation family we will show that every compact parallelizable complex manifold $X$ with $b_1(X) > 0$ admits non-trivial small deformations.

**Theorem 14.** *Let $G$ be a connected complex Lie group, $\Gamma$ a discrete cocompact subgroup and $X = G/\Gamma$ the quotient manifold. Assume that $b_1(X) > 0$.*

*Then there exist small deformations of $X$, i.e. there is a proper flat holomorphic family $Y \to \Delta_1$ with $Y_0 \simeq X$ and $Y_t \not\simeq X$ for $t \neq 0$.*

*Proof.* We will reduce the general theorem to two special cases, namely $G \simeq SL_2(\mathbb{C})$ and the case where $G$ is commutative.

**Claim 5.** *Under the above assumptions there exists a connected complex Lie group $H$ which is either commutative or isomorphic to $SL_2(\mathbb{C})$, a cocompact discrete subgroup $\Lambda \subset H$ and a $G$-equivariant holomorphic surjection $\pi : X \to Z = H/\Lambda$. Furthermore $b_1(Z) > 0$.*

*Proof.* Let $R$ denote the radical of $G$. By a theorem of Mostow there is a fibration $\rho : G/\Gamma \to G/R\Gamma$. If $b_1(G/R\Gamma) > 0$, arithmeticity results for lattices in semisimple Lie groups (see [14]) imply that there exists a fibration $\pi_0 : G/R\Gamma \to SL_2(\mathbb{C})/\Lambda$ such that $\pi = \pi_0 \circ \rho$ is the desired surjection. On the other hand, if $b_1(G/R\Gamma) = 0$, then $G'\Gamma$ is closed in $G$ and $\Gamma'$ is a subgroup of finite index in $G' \cap \Gamma$. Thus in this case $G/\Gamma \to G/G'\Gamma$ yields a surjection onto a positive-dimensional torus. □

**Claim 6.** *Let $G$ be either $SL_2(\mathbb{C})$ or $(\mathbb{C}^n, +)$, $\Gamma$ a discrete cocompact subgroup, $X = G/\Gamma$ the quotient manifold and $\tau : \Gamma \to (\mathbb{Z}, +)$ a surjective group homomorphism.*

*Then there exists an open neighbourhood $U$ of $e$ in $G$ such that for every $u \in U$ the $\Gamma$-action on $G$ given by*

$$\gamma : x \mapsto u^{-\tau(\gamma)} \cdot x \cdot \gamma$$

*is free and properly discontinuous. Moreover $U$ contains a subset $A$ of measure zero such that the quotient manifold $X_u$ is not biholomorphic to $X$ for any $u \notin A$.*

*Proof.* For $G = (\mathbb{C}^n, +)$ this is easy to check and for $G = SL_2(\mathbb{C})$ it has been proved by Ghys [7]. □

Thus we deduced that we have a fibration $G/\Gamma \to H/\Lambda$ and non-trivial small deformations of $H/\Lambda$. We need to show that these deformations can be lifted to non-trivial deformations of $G/\Gamma$.

This is achieved by the following statement.

**Claim 7.** *Let $G$, $H$ be connected Lie groups, $\Gamma$ resp. $\Lambda$ discrete cocompact subgroups in $G$ resp. $H$, $\pi : G \to H$ a surjective Lie group homomorphism with $\pi(\Gamma) = \Lambda$, $\tau : \Lambda \to \mathbb{Z}$ a surjective group homomorphism and $u \in G$.*

*Assume that the $\Lambda$-action on $H$ given by*

$$\lambda : h \mapsto (\pi(u))^{-\tau(\lambda)} \cdot h \cdot \lambda$$

*is free and properly discontinuous.*

*Then the $\Gamma$-action on $G$ given by*

$$\gamma : u^{-\tau(\pi(\gamma))} \cdot g \cdot \gamma$$

*is also free and properly discontinuous.*

*Proof.* In order to show that the action on $G$ is properly discontinuous, we have to verify that for every compact subset $K \subset G$ the set

$$S = \{\gamma \in \Gamma : u^{-\tau(\pi(\gamma))} \cdot K \cdot \gamma \cap K \neq \emptyset\}$$

is finite. Now $\pi(K)$ is compact and hence
$$S_1 = \{\lambda \in \Lambda : (\pi(u))^{-\tau(\lambda)} \cdot \pi(K) \cdot \lambda \cap \pi(K) \neq \emptyset\}$$
must be finite. Thus $M = \{u^{-\tau(\lambda)} : \lambda \in S_1$ is likewise finite. Now
$$S \subset \{\gamma \in \Gamma : M \cdot K \cdot \gamma \cap K \neq \emptyset\}.$$
Both $M \cdot K$ and $K$ are compact and consequently $S$ is finite. The freeness can be checked in a similar way. □

□

## 9. Examples

1. Let $\Lambda \subset SL_2(\mathbb{C})$ be a cocompact lattice with a surjective group homomorphism $\rho : \Lambda \to \mathbb{Z}$ (see e.g. [12] for existence of such lattices). Let $G = SL_2(\mathbb{C}) \times \mathbb{C}$ and $\Gamma = \{(s,x) \in G : s \in \Gamma, x - \sqrt{2}\rho(s) \in \mathbb{Z}[i]\}$. Then $\Gamma$ is a lattice in $G$ such that $G'\Gamma$ is not closed in $G$. Furthermore $b_1(G/\Gamma) = rank_\mathbb{Z}(\Lambda/\Lambda') + 2 > 0$ and $\dim H^1(G/\Gamma, \mathcal{O}) = rank_\mathbb{Z}(\Lambda/\Lambda') + 1 > 0$ although every holomorphic map from $G/\Gamma$ to a torus is constant.

2. Nakamura ([16]) gave an example of a three-dimensional solvmanifold $X = G/\Gamma$ with $\dim G/G' = 1$ and $\dim H^1(G/\Gamma, \mathcal{O}) = 3$. We will give another such example which in addition demonstrates that $\dim H^1(G/\Gamma, \mathcal{O})$ may jump within a commensurability class of $\Gamma$.

Let $p$ be a non-square positive natural number, $L = \mathbb{Q}[\sqrt{p}]$ and $K = L[i]$. By Dirichlet's theorem the group $\mathcal{O}_L^*$ of units of $L$ contains an element of infinite order $\alpha$. We may assume that $N_{L/\mathbb{Q}}(\alpha) = 1$. Fix two embeddings $\sigma_1, \sigma_2 : K \to \mathbb{C}$ such that $\sigma_1 \neq \overline{\sigma_2}$, but $\sigma_1(i) = i$ and $\sigma_2(i) = -i$. Now $\sigma = (\sigma_1, \sigma_2) : K \to \mathbb{C}^2$ embedds $\mathcal{O}_K$ in $\mathbb{C}^2$ as a lattice. Let $\Lambda_i$ $(i = 1,2)$ be the subgroups of $\mathcal{O}_K^*$ generated by $\alpha$ resp. $\alpha$ and $i$. Now $\sigma$ induces an embedding of the groups $\Lambda_i$ into $GL(2,\mathbb{C})$. Note that $\det \sigma(\alpha) = N_{L/\mathbb{Q}}(\alpha) = 1$ and $\det \sigma(i) = 1$. Hence both $\sigma_i(\Lambda_i)$ are lattices in a Cartan subgroup $T$ of $SL_2(\mathbb{C})$. Let $\Gamma_i = \Lambda \ltimes \mathcal{O}_K$. Then the groups $\Gamma_i$ are lattices in a three-dimensional solvable complex Lie group $G = \mathbb{C}^* \ltimes (\mathbb{C}^2, +)$. By construction, the $\Gamma_1$-action on $\mathbb{C}^2$ is totally real, while $\sigma(i) \in \Gamma_2$ acts on $\mathbb{C}^2$ without real eigenvalues.

Therefore we obtain a two-to-one covering of three-dimensional solvmanifolds
$$G/\Gamma_1 \to G/\Gamma_2$$
with $\dim H^1(G/\Gamma_1, \mathcal{O}) = 3$ and $\dim H^1(G/\Gamma_2, \mathcal{O}) = 1$.

JÖRG WINKELMANN, MATH. INSTITUT NA 4/69, RUHR-UNIVERSITÄT BOCHUM, 44780 BOCHUM, GERMANY
*E-mail address*: jw@cplx.ruhr-uni-bochum.de